\documentclass[10pt]{amsart}
\usepackage{amsmath,mathrsfs}
\usepackage{amssymb}
\usepackage{latexsym}
\usepackage{amsfonts}
\usepackage{graphicx}
\usepackage{psfrag}

\newtheorem{Pa}{Paper}[section]
\newtheorem{theorem}[Pa]{{\bf Theorem}}
\newtheorem{lemma}[Pa]{{\bf Lemma}}
\newtheorem{corollary}[Pa]{{\bf Corollary}}

\newtheorem{proposition}[Pa]{{\bf Proposition}}
\newtheorem{observation}[Pa]{{\bf Observation}}

\newtheorem{example}[Pa]{{\bf Example}}

\def\C{\mathbb C}
\def\R{{\mathbb R}}

\def\@tempb{saamsart}
\usepackage{hyperref}
\begin{document}
\bibliographystyle{plain}
\thispagestyle{empty}
\title[
realization of generalized
positive rational functions]
{The positive real lemma and construction
of all realizations of generalized positive rational
functions}
\author[D. Alpay]{Daniel Alpay}
\address{(DA) Department of mathematics,
Ben-Gurion University of the Negev, P.O. Box 653, Beer-Sheva
84105, Israel} \email{dany@math.bgu.ac.il}
\author[I. Lewkowicz]{Izchak ~Lewkowicz}
\address{(IL) Department of electrical
engineering, Ben-Gurion University of the Negev, P.O. Box 653,
Beer-Sheva 84105, Israel} \email{izchak@ee.bgu.ac.il}

\thanks{D. Alpay thanks the
Earl Katz family for endowing the chair which supported his
research. This research is part of the European Science
Foundation Networking Program HCAA, and was supported in part by
the Israel Science Foundation grant 1023/07}


\def\squarebox#1{\hbox to #1{\hfill\vbox to #1{\vfill}}}

\maketitle
\begin{abstract}
We here extend the well known Positive Real Lemma (also known as
the Kalman-Yakubovich-Popov Lemma) to complex matrix-valued
generalized positive rational function, when non-minimal
realizations are considered. All state space realizations are
partitioned into subsets, each is identified with a set of
matrices satisfying the same Lyapunov inclusion. Thus, each
subset forms a convex invertible cone, {\bf cic} in short, and
is in fact is replica of all realizations of positive functions
of the same dimensions. We then exploit this result to provide
an easy construction procedure of all (not necessarily
minimal) state space realizations of generalized positive
functions. As a by-product, this approach enables us to
characterize systems which can be brought, through static output
feedback, to be generalized positive.
\end{abstract}

\section{Introduction}
\setcounter{equation}{0}
 For a half of a century, the Positive Real Lemma (also known as
the Kalman-Yakubovich-Popov Lemma) has been recognized as a
fundamental result in System Theory. We here extend and exploit
it in various ways. Let $\C_+$ and $\C_-$ be the open right and
left halves of the complex plane respectively, and $~{\mathbb
P}_k,~ (\overline{\mathbb P}_k)~$ be the set of all $~k\times k~$
positive definite (semidefinite) matrices\begin{footnote}
{Whenever clear from the context, the dimension subscript will be
omitted.}\end{footnote}. Recall that a $~p\times p$-valued
function $~F(s)$, analytic in $\mathbb C_+$ is said to be {\em
positive} if
\begin{equation}\label{eq:PosFunc}
\begin{matrix}
F(s)+F(s)^*
\in\overline{\mathbb P}_p &~&s\in\C_+~.\end{matrix}
\end{equation}
The study of rational positive functions, denoted by
$~{\mathcal P}$, has been motivated from the 1920's by (lumped)
electrical networks theory, see e.g. \cite{AV}, \cite{Be}. From
the 1960's positive functions also appeared in books on absolute
stability theory, see e.g. \cite{NT}, \cite{Po}.

A $~p\times p$-valued function of bounded type in $\mathbb C_+$
(i.e. a quotient of two functions analytic and bounded in
$\mathbb C_+$) is called {\em generalized positive}
$~\mathcal{GP}~$ if
\begin{equation}\label{DefGP}
F(i\omega)+F(i\omega)^*
\in\overline{\mathbb P}_p
\quad a.e.\quad\quad\omega\in\R,
\end{equation}
where $F(i\omega)$ denotes the non-tangential
limit\begin{footnote}
{This limit exists almost everywhere on $i\R$ because $F$ is
assumed of bounded type in $\mathbb C_+$, see e.g.
\cite{Du}.}
\end{footnote}
of $F$ at the point $i\omega$.

Generalized positive functions were introduced in the context
of the Positive Real Lemma (PRL), see \cite{AM} and references
therein\begin{footnote}{The original formulation was real.
The case we address is in fact {\em generalized} positive and
{\em complex}, but we wish to adhere to the commonly used
term: Positive Real Lemma.}\end{footnote}. Applications of
$~\mathcal{GP}~$ functions to electrical networks appeared in
\cite{IO}, and to control in \cite{LSC}, where they first
casted in a Linear Matrix Inequality (LMI) framework, see e.g.
\cite{BGFB} for more information on LMI. For more application
of the generalized PRL, see \cite{HSK}
\vskip 0.2cm

Both function sets $\mathcal P$ and $\mathcal {GP}$ are closed under
positive scaling, sum and inversion (when the given function has
a non-identically vanishing determinant). Yet these spaces have
quite different properties, as we now illustrate. Let $\mathbb
C^{p\times p}(X)$ be the space of $\mathbb C^{p\times p}$-valued
rational functions. Recall that a function $\Psi\in \mathbb
C^{p\times p}(X)$ belongs to $\mathcal {GP}$ if and only if it
can be factorized as
\begin{equation}
\label{eq:factor}
\Psi(s)=G(s) P(s) G(-s^*)^*,
\end{equation}
where $~G,P \in\mathbb C^{p\times p}(X)$, and
$~P\in\mathcal P$. Factorization of this nature appeared e.g.
in \cite{DHdS} and \cite{DLLS1}, see also Observation
\ref{GP-L(I)} below. The significance of \eqref{eq:factor} to
scalar rational $\mathcal{GP}$ functions was recently treated
in \cite{AL1} and \cite{AL2}.
\vskip 0.2cm

We now consider properties of the sum of two rational
$\mathcal{GP}$ functions (series connection in electrical
engineering jargon). From \eqref{DefGP} it follows that this
sum is again in $\mathcal{GP}$, but both the McMillan degree
and the number of negative squares (roughly, the number of
poles in $\C_+$) increase. For more details on the number of
negative squares see \cite{krein} and \cite{KL1}. Recall that
a rational function $\Psi$ is in $\mathcal{GP}$ if and only if
the kernel
\[
\frac{\Psi(s)+\Psi(w)^*}{s+w^*}
\]
has a finite number of negative squares in its domain of
definition in $\mathbb C_+$. The number of negative squares of
the sum of two elements in $\mathcal {GP}$ is preserved if, for
instance, in \eqref{eq:factor} $\Psi_1(s)=G(s) P_1(s) G(-s^*)^*$
and $\Psi_2(s)=G(s)P_2(s) G(-s^*)^*$ with $P_1, P_2\in{\mathcal P}$
and the same function $G\in\mathbb C^{p\times p}(X)$, see
\cite[Section 3]{AL2} for the scalar case.

In contrast, if one takes a {\em state space realization sum}, the
McMillan degree of the resulting function does not increase, but
it may turn to not generalized positive at all, see Example
\ref{ConvRealiz} below. One of the results of this paper
is a partitioning of $\mathcal{GP}$ functions to subsets, denoted
by $\mathcal{GP}(r, \nu, p)$, closed under state space addition,
while both the McMillan degree and the number of negative squares
do not increase. See Theorem \ref{GP(r,nu,p)cic}
below.
\vskip 0.2cm

We resort to some preliminaries. Let
$\Psi\in\mathbb C^{p\times p}(X)$ be
of McMillan degree $~q$, and analytic at infinity, i.e.
$~\lim\limits_{s~\rightarrow~\infty}~\Psi(s)~$ exists. Namely,
$\Psi$ admits a state space realization
\begin{equation}\label{StateSpace}
\begin{matrix}
\Psi(s)=C(sI-A)^{-1}B+D&~&~&
L:=\begin{pmatrix}A&~B\\C&~D
\end{pmatrix}\end{matrix}
\end{equation}
with $~A\in\C^{n\times n}$, $n\geq q$, $~B, C^*\in\C^{n\times
p}~$ and $~D\in\C^{p\times p}$, namely, $~L\in\C^{(n+p)
\times(n+p)}$. If the McMillan degree of $\Psi(s)~$ satisfies
$~q=n$, the realization is called minimal.

We can now state the Positive Real Lemma (PRL) as
presented in \cite[Theorem 1]{DDGK}
(up to substituting the real setting by a complex one)
\vskip 0.2cm

\begin{lemma}\label{gpLemma}
Let $~\Psi(s)~$ be a $~p\times p$-valued rational function in
\eqref{StateSpace} and assume that $~q=n$.

(I) $~\Psi\in\mathcal{GP}~$ if and only if
\begin{equation}\label{Lyap}
HL+L^*H=Q\in\overline{\mathbb P}_{n+p}~,
\end{equation}
for some $~H={\rm diag}\{\hat{H},~I_p\}$, where
$~\hat{H}\in\C^{n\times n}~$ is Hermitian
non-singular.

(II) If $~\Psi\in{\mathcal P}~$ then in part (I)
$~-\hat{H}\in{\mathbb P}_n~$.
\end{lemma}
\vskip 0.2cm

The aim of this work is to first extend this result to the
non-minimal case, and then to use it to obtain a straightforward
construction of ~{\em all}~ (not necessarily minimal) state space
realization of $~{\mathcal P}~$ and $~\mathcal{GP}~$ rational
functions. This is then used to describe the already mentioned
partitioning to sets of the form $~\mathcal{GP}(r, \nu, p)$ and
to obtain other results, described below.
\vskip 0.2cm

The outline of the paper is as follows: The paper is composed of
eight sections besides the introduction. In Section \ref{Sec:2}
we give a short review of the literature, which should be
paralleled with a complementary survey we offered in our previous
paper \cite{AL1}. Our aim is not to provide a complete survey,
but to raise, through samples, the intriguing observation that
although the PRL has been a standard textbook material from the
1970's, see e.g. \cite{AV}, \cite[Chapter 3]{FCG},
\cite[Section 4.4, Appendix]{NT} and \cite[Section 8.5]{Po}, it
is not straightforward to cover the relevant literature. In simple
words, there is too little of cross-referencing
\begin{footnote}{For example, it seems that \cite{DDGK}
was hardly ever cited.}\end{footnote}.

In Section \ref{Sec:3} an algebraic
Riccati inclusion associated with necessity part of the
generalized positive real Lemma is addressed. Sections
\ref{Sec:4} and \ref{Sec:5} are devoted to showing that an
algebraic Lyapunov inclusion associated with the sufficiency part
of the generalized positive real Lemma is independent of the
minimality of the realization. Some background material
concerning sets of Lyapunov inclusions is reviewed in Section
\ref{Sec:6}. In particular we provide a convenient parameterization
of all matrices $L$ satisfying the Lyapunov inclusion \eqref{Lyap}
where $H$ is fixed and $Q$ varies over all $\overline{\mathbb P}$.
This Lyapunov inclusion formulation is then employed in Section
\ref{Sec:7} to provide a straightforward parameterization of all
state space realization of $~\mathcal{GP}~$ rational functions.
This allows us to describe $~\mathcal{GP}~$ functions as a
union of replicas of positive functions.
As an application of theses sets, in Section \ref{Sec:8} we
characterize all rational functions (vanishing at infinity)
which can be rendered $~\mathcal{GP}$, through a static state
feedback.
Concluding remarks are given in Section \ref{Sec:10}.

\section{A historical perspective}
\label{Sec:2}
\setcounter{equation}{0}
We here review some of the
relevant existing literature.

As mentioned, the above version of Lemma \ref{gpLemma} is
from \cite{DDGK} and was repeated in \cite{LSC}. It was
originally proved in \cite{AM}. A special case was later
treated in \cite{MSS}.  The positive function case
(part II) is well known and sometimes is referred to as the
Kalman-Yakubovich-Popov Lemma and is dated to the 1960's.
For an early full account see e.g. \cite[Chapter 5]{AV}. An
easy-to-read historical perspective is given in
\cite[Sections 1]{BGFB}.

A matrix formulation, see \eqref{StateSpace}, of the PRL
was introduced in the PhD. thesis of P. Faurre, see
\cite[Theorem 4.2]{F} and then in a book he co-authored,
\cite[Theorem 3.1]{FCG}. In fact it implicitly earlier
appeared in \cite{W}. The formulation through the
Rosenbrock system matrix $~L~$ \eqref{StateSpace} (for
not necessarily positive systems) explicitly introduced
in \cite{DDGK} and subsequently in \cite[Lemma 8]{LSC}.
An interesting special case was studied in
\cite[Theorem 4]{GvDKDM}.
The notion of Linear Matrix Inequality (LMI) was introduced
in \cite{W}. \cite{LSC} was one of the early works recognizing
the applicability of LMI framework to the (generalized)
Positive Real Lemma (PRL), see also \cite[Section 2.7.2]{BGFB}.
A comprehensive survey of the LMI
approach to the PRL appeared in \cite{GL}. Unfortunately, in
spite of its admirable reference list (201 items), some important
relevant results are missing.
\vskip 0.2cm

Following \eqref{eq:PosFunc} Positive functions map $~\C_+~$
to $~\overline{\mathbb L}(I)$, the set of matrices with a
non-negative Hermitian part\begin{footnote}{The notation
$\overline{\mathbb L}(I)$, with ${\mathbb L}$ honoring A.M.
Lyapunov, will be formally defined in \eqref{CommonLyap}
below.}\end{footnote}. Analogously,
following \eqref{DefGP}, a $~\mathcal{GP}~$ function maps
$~i\R~$ (after removing all poles of the function) into
$\overline{\mathbb L}(I)$, see Observation \ref{GP-L(I)} below.
Closely related function sets are addressed
in the literature:

\begin{itemize}
\item{}Bounded functions mapping $~\C_+~$ to weak
contractions while generalized bounded functions map
$~i\R~$ to weak contractions. Versions of the PRL
for bounded functions appeared \cite[Theorem 4.2]{F}
\cite[Section 7.2]{AV} and for generalized bounded
in \cite[Equation (6)]{DDGK}. An interesting subclass
is treated in \cite[Theorem 2.1]{AG1}.
An important result on a subclass of
generalized bounded functions where $~i\R~$ is
substituted by $~\R$, appeared in \cite[Theorem 3.2]{GoRu}.
See also \cite[Theorem 2.12]{AG2} and \cite[Chapter 21]{LR}.

\item{}Carath\'eodory functions, analytically mapping
the open unit disk to $\overline{\mathbb L}(I)$ and
generalized (=pseudo) Carath\'{e}odory functions,
mapping the unit circle (from the poles of the given
function have been removed) to $\overline{\mathbb L}(I)$.
Versions of the PRL for Carath\'eodory functions appeared
in \cite{XH}.

\item{}Schur functions, analytically mapping the open
unit disk to weak contractions and generalized Schur
functions, mapping the unit circle to weak contractions.
Version of the PRL for Schur functions appeared in
\cite[Theorem 2.5]{AG3}, \cite[Section II]{V}, \cite{XH}
and for operator valued in \cite{Ar}. For generalized
Schur functions see \cite[Theorem 2]{DDGK} and
for an interesting special case see \cite[Theorem 3]{GvDKDM}.
A time-varying extension of the PRL to generalized Schur
functions along with a thorough study of various applications
is are given in \cite[Theorem 1.2.5 and Appendix 3A]{HSK}.

\item{}Nevanlinna functions, analytically mapping the open
upper half plane to $~i\overline{\mathbb L}(I)$ and the
generalized Nevanlinna functions map $~\R~$ to
$~i\overline{\mathbb L}(I)$, see e.g. \cite{DHdS} and
\cite{DLLS1}. Minimal realization of infinite dimensional
generalized Nevanlinna functions was studied in \cite{DLLS2}.
\end{itemize}

As already stated, we do not aspire to provide a survey of PRL
related results, and we are aware of additional references dealing
with the subject, not mentioned here. We focused on the scattered
nature of the literature related to the ~{\em generalized}~ case.\\

\section{Generalized positive lemma necessity and
the Riccati equation}
\setcounter{equation}{0}
\label{Sec:3}
It has been long recognized that with the part (b) of Lemma
\ref{gpLemma} (dealing with positive functions) one can associate
an algebraic Riccati equation, see e.g. \cite[Section 5.4]{AV}
\cite[Chapter 5]{FCG} and \cite[Section 2.7.2]{BGFB}. For Schur
functions version see e.g. \cite[Theorem 2.1]{deSX}. We now
address the extension of this result to $~\mathcal{GP}~$
functions. It first appeared in the context of a variant of
Generalized Bounded functions (where $~i\R~$ was substituted by
$~\R$) in \cite[Theorem 3.2]{GoRu} and in \cite[Theorem
2.12]{AG2}. For another variant of this result, see \cite[Section
20.1]{LR}. In Proposition \ref{gpLemmaNecessity1} below we
provide a simple proof of the result using the system matrix
formulation \eqref{StateSpace}, employed all along this work. We
shall find it convenient to resort to the following notation of
sets of matrices sharing a common Lyapunov factor: For a
$~r\times r~$ Hermitian non-singular matrix $~H$, define the sets
of $~r\times r~$ matrices, $~{\mathbb L}(H)~$ and
$~\overline{\mathbb L}(H)$ as,
\begin{equation}\label{CommonLyap}
\begin{matrix}
(a)~{\mathbb L}(H):=\{ L~:~HL+L^*H\in{\mathbb P}_r~
\}
&~&
(b)~\overline{\mathbb L}(H):=\{ L~:~
HL+L^*H\in\overline{\mathbb P}_r~\}.
\end{matrix}
\end{equation}
In particular, $(\overline{\mathbb L}(I)~)$ $~{\mathbb L}(I)~$
is the set of matrices with positive (semi)definite Hermitian part.

\begin{proposition}\label{inertia4}
Consider \eqref{StateSpace}, \eqref{Lyap} with
$~H={\rm diag}\{\hat{H},~I_p\}$ and $~\hat{H}\in\C^{n\times n}$
Hermitian. Assume in addition that $~D\in{\mathbb L}(I_p)$.
Let us define the following $n\times n$ Riccati expression,
\begin{equation}\label{Riccati}
\begin{matrix}
M&:=\hat{H}(A-B(D+D^*)^{-1}C)&+&(A-B(D+D^*)^{-1}C)^*\hat{H}\\~
&-\hat{H}B(D+D^*)^{-1}B^*\hat{H}&-&C^*(D+D^*)^{-1}C.
\end{matrix}
\end{equation}
Then $~Q\in\overline{\mathbb P}_{n+p}$, if
and only if in \eqref{Riccati}
\begin{equation}\label{RiccatiInc}
M\in\overline{\mathbb P}_n~.
\end{equation}
\end{proposition}

{\bf Proof}\quad
Using the fact that $~D\in{\mathbb L}(I_p)~$ one can employ
the classical Schur's complement, e.g. \cite[Theorem 7.7.6]{HJ1},
to write down $~Q~$ in \eqref{Lyap} explicitly,
$$
Q=\begin{pmatrix}
\hat{H}A+A^*\hat{H}&\hat{H}B+C^*\\C+B^*\hat{H}&D+D^*
\end{pmatrix}
=\begin{pmatrix}I_n&R
\\ 0&I_p\end{pmatrix}
\begin{pmatrix}M&0\\ 0&D+D^*\end{pmatrix}
\begin{pmatrix}I_n&0\\ R^*
&I_p\end{pmatrix},
$$
where $~R:=(\hat{H}B+C^*)(D+D^*)^{-1}$ and $~M~$ is given
in \eqref{Riccati}. Thus indeed
$~Q\in\overline{\mathbb P}_{n+p}$
if and only if $~M\in\overline{\mathbb P}_n$.
\qed
\vskip 0.2cm

We can now re-formulate the necessity part of Lemma \ref{gpLemma}.

\begin{proposition}\label{gpLemmaNecessity1}
Let $~\Psi\in\mathcal{GP}~$ be a $~p\times p$-valued
rational function so that
$~\lim\limits_{s~\rightarrow~\infty}~\Psi(s):=D~$ exists.
Furthermore assume that $~D\in{\mathbb L}(I_p)$. Let $~q~$ be the
McMillan degree of $~\Psi$.

(I) Then, $~\Psi(s)~$ admits a state space realization
\eqref{StateSpace}, so that it is minimal $(q=n)$ and it satisfies
the Riccati inclusion in \eqref{Riccati}, \eqref{RiccatiInc} for
some $~n\times n~$ Hermitian non-singular \mbox{$\hat{H}$.}

(II) If $~\Psi\in{\mathcal P}~$ then in part (I)
$~-\hat{H}\in{\mathbb P}_n$.
\end{proposition}

Note that the technical condition $~D\in{\mathbb L}(I_p)~$
in Proposition \ref{gpLemmaNecessity1} is indeed
restrictive.
For example, many system of interest have a zero at
infinity and thus are excluded from the discussion. On the
other hand, whenever \mbox{$\Psi(s)=C(sI-A)^{-1}B+D$} is a
$~{\mathcal G}{\mathcal P}~$ function, from the above
discussion it follows that $~D\in\overline{\mathbb L}(I_p)$,
hence one can always construct another $~\mathcal{GP}~$
function \mbox{$\tilde{\Psi}(s)=C(sI-A)^{-1}B+\tilde{D}$,}
with \mbox{$\tilde{D}\in{\mathbb L}(I_p)$} so that
$~\epsilon\geq\| D-\tilde{D}\|$, where $~\epsilon>0~$ is
arbitrary.

\section{Generalized positive lemma sufficiency - an
extension}\setcounter{equation}{0}
\label{Sec:4}

The sufficiency statement of the (generalized) Positive
Real Lemma (PRL) of matrix valued rational functions
was first proved in \cite{AM}, under the assumption
of minimality of the realization ($q=n$). We now address
the question of relaxing this minimality constraint.
This problem was treated in the framework of
positive functions in \cite[Lemma 6]{XH} and in
the framework of functions satisfying
$~D\in{\mathbb L}(I_p)~$ (as in the previous section) in
\cite[Section 21.3]{LR}. In a different formulation see also
\cite[Theorem 1]{Ra}. In \cite{FP} a proof of the
sufficiency statement of the (generalized) PRL, removing
the minimality of realization condition, was presented.
Unfortunately a (redundant) spectral condition on $~A~$
was imposed there. The result of Proposition
\ref{gpLemmaSufficiency} below, avoids any restriction.

In addition, in Proposition \ref{gpLemmaSufficiency}
below we show that one can bound the number of poles of
$~\Psi(s)~$ in each open half plane\begin{footnote}{
Although in a different framework,
bounds of a similar nature can be found in
\cite[Theorem 3.4]{GoRu} and
subsequently in \cite[Theorem 2.12]{AG2}
and in \cite[Section 21.2]{LR}.}\end{footnote}.
To this end, we need some preliminaries.
Recall that for a matrix $~A\in\C^{n\times n}$ one
can associate a triple:
\mbox{${\rm inertia}(A)=(\nu,~\delta,~\pi)$,}
with $~\nu+\delta+\pi=n$, if
$~A~$ has $~\nu~$ eigenvalues in $~\C_-$,
$~\pi~$ eigenvalues in $~\C_+~$ and $~\delta~$ eigenvalues
on $~i\R$, see e.g. \cite[2.1.1]{HJ2}.
Let $~A,\hat{H}\in\C^{n\times n}$ with $~\hat{H}~$ Hermitian,
be with inertia,
\begin{equation}\label{eq:inertia}
{\rm inertia}(A)=(\nu_A,~\delta_A,~\pi_A)
\quad\quad
{\rm inertia}(\hat{H})=(\nu,~0,~n-\nu)
\quad\quad \nu\in[0,~n],
\end{equation}
i.e. $\hat{H}$ is non-singular.
Consider now the Lyapunov equation
\begin{equation}\label{Lyap1}
\hat{H}A+A^*\hat{H}=\hat{Q}\in\overline{\mathbb P}_n~.
\end{equation}
Recall that from a pair $~A\in\C^{n\times n}$,
$B\in\C^{n\times p}$, the following controllability
matrix may be constructed
\mbox{${\mathcal C}:=[B~\vdots~AB~\vdots~\cdots~
\vdots~A^{n-1}B]$}. Then $~X_{\rm cont}(A, B)$, the
controllable subspace associated with the pair
$~A, B~$ is given by the range of $~{\mathcal C}$ and
$~X_{\rm cont}(A, B)^{\perp}$, its orthogonal complement,
is given by the null-space of $~{\mathcal C}^*$,
see e.g. \cite[Definition 2.4.8]{HJ2}. Similarly
with a pair $~A\in\C^{n\times n}$, $C\in\C^{p\times n}$,
one can associate a observable subspace
$~X_{\rm obs}(A, C)$, given by
$~X_{\rm obs}(A, C)=X_{\rm cont}(A^*, C^*)$ (and
$~X_{\rm obs}(A, C)^{\perp}=X_{\rm cont}(A^*, C^*)^{\perp}$).
For $~A,~\hat{Q}~$ in \eqref{Lyap1} denote,
\begin{equation}\label{eq:m}
m:={\rm dim}~X_{\rm obs}(A,~\hat{Q})^{\perp}.
\end{equation}
Namely, pair $~A,~\hat{Q}~$ is observable whenever $~m=0$.
We can now cite the following important result
of R. Loewy \cite{Lo}, adapted to our framework,

\begin{theorem}\label{Th:Loewy}
Let $~A, \hat{H}~$ and $~m~$ be as in \eqref{eq:inertia},
\eqref{Lyap1} and \eqref{eq:m}. Then,
$$\begin{matrix}
\nu\geq\nu_A\geq\max(0,~\nu-m)&~&~&n-\nu\geq\pi_A\geq
\max(0,~n-\nu-m).
\end{matrix}$$
\end{theorem}
We can now state the extended sufficiency part of
the PRL.

\begin{proposition}\label{gpLemmaSufficiency}

Let $~\Psi~$ in \eqref{StateSpace} be a $~p\times p$-valued
rational function of McMillan degree $~q$, $n\geq q$.
Assume that state space realization in \eqref{StateSpace}
satisfies the Lyapunov equation \eqref{Lyap} with
$~H={\rm diag}\{\hat{H},~I_p\}$, where $~\hat{H}~$ is
$~n\times n~$ Hermitian.

(I) If $~{\rm inertia}(\hat{H})=(\nu,~0,~n-\nu)$, for
some $~\nu\in[0,~n]$, then, $~\Psi~$ is a $~\mathcal{GP}~$
function with at most $~\nu~$ poles in $~\C_-~$ and
$~n-\nu~$ poles in $~\C_+$.

(II) If in part (I) $~-\hat{H}\in{\mathbb P}_n$, i.e. $~\nu=n$,
then $~\Psi\in{\mathcal P}$.
\end{proposition}

{\bf Proof :}\quad I. Indeed assume that \eqref{Lyap} is
satisfied with $~H={\rm diag}\{\hat{H}, I_p\}$, $~\hat{H}~$
Hermitian nonsingular and $~L~$ as in \eqref{StateSpace}. Note
that in \eqref{Lyap} $~Q~$ is in $\overline{\mathbb P}_{n+p}$,
thus its upper left block is in $\overline{\mathbb P}_n$.
Namely, \eqref{Lyap1} is satisfied, so by Theorem \ref{Th:Loewy}
the matrix $~A~$ has at most $~\nu~$ eigenvalues in $~\C_-~$ and
$~n-\nu~$ eigenvalues in $~\C_+$. Recall that the poles of
$~\Psi~$ are determined by the eigenvalues of $~A$, see
\eqref{StateSpace}.

Next denoting $~S:={\rm diag}\{-sI_n, 0_p\}~$ we note that
\begin{equation}\label{LyapS}
HS+S^*H
\in\overline{\mathbb P}_{n+p}~,
\end{equation}
for all $~s\in{i\R}$. Take now \begin{footnote}{It is
interesting to note that the zeroes of $~\Psi(s)~$ are
the points $~s~$ for which $~\tilde{L}~$ is singular.}
\end{footnote}
$~\tilde{L}:=L+S=${\mbox{
\tiny$\begin{pmatrix}-sI_n+A&~B\\C&~D\end{pmatrix}$}}.
Then for all $~s\in{i\R}$ also,
\begin{equation}\label{LyapTildeL}
H\tilde{L}+\tilde{L}^*H
\in\overline{\mathbb P}_{n+p}~.
\end{equation}
Next, recall that for arbitrary constant matrix
$~\Psi\in\mathcal{GP}$, if and only if
\mbox{$(\Psi+T)\in{\mathcal G}{\mathcal P}$,} for
arbitrary constant
matrix $~-T^*=T\in\C^{p\times p}$. Thus, up to a shift by a
skew-Hermitian matrix, we can assume that $~\Psi(s)~$ in
\eqref{StateSpace} is almost everywhere invertible. Thus, a
straightforward calculation (see e.g.
\cite[0.7.3]{HJ1}) results in,
$$\tilde{L}^{-1}=
\begin{pmatrix}
(sI-A)^{-1}\left(B\Psi(s)^{-1}C(sI-A)^{-1}-I\right)~~&
(sI-A)^{-1}B\Psi(s)^{-1}\\ \Psi(s)^{-1}C(sI-A)^{-1}&
\Psi(s)^{-1}\end{pmatrix}.$$
Multiplying \eqref{LyapTildeL} by
$~\left(\tilde{L}^*\right)^{-1}~$ from the left and
$~\tilde{L}^{-1}~$ from the right, yields
\begin{equation}\label{LyapRatio}
H\tilde{L}^{-1}+(\tilde{L}^{-1})^*H
\in\overline{\mathbb P}_{n+p}
\end{equation}
for all $~s\in{i}\R$. Now in particular the $~p\times p~$
lower right block of \eqref{LyapRatio} satisfies,
\begin{equation}\label{defGP}
\Psi(s)^{-1}+\left(\Psi(s)^{-1}\right)^*=
\left(
H\tilde{L}^{-1}+(\tilde{L}^{-1})^*H\right)_{22}\in
\overline{\mathbb P}_p
\end{equation}
for all $~s\in{i}\R$. Thus, $~\Psi(s)^{-1}~$ is in
$~{\mathcal G}{\mathcal P}~$ and hence also $~\Psi(s)$.
Thus, the first part of the claim is established.

(b) Positive functions. If $~-\hat{H}\in\overline{\mathbb P}_n~$
the relation in \eqref{LyapS} holds for all $~s\in\C_+~$
and subsequently, also \eqref{LyapRatio} and \eqref{defGP}.
Hence, $~\Psi\in{\mathcal P}$, so the proof is complete.
\qed
\vskip 0.2cm

In the next section we scrutinize some aspects of
Proposition \ref{gpLemmaSufficiency}.

\section{non-minimal realization and bounds on
inertia - a closer look}
\setcounter{equation}{0}
\label{Sec:5}

Roughly, application of Theorem \ref{Th:Loewy} to Proposition
\ref{gpLemmaSufficiency} suggests that the ``further"
from minimality the realization is, the cruder is the bound
on the number of poles in $\C_+$. We here illustrate the
``at most" clause in the
statement Proposition \ref{gpLemmaSufficiency} with respect
to the number of poles in each open half plane.

\begin{example}
\label{InertiaExample}
{\rm All functions considered in this example are so that in
\eqref{Lyap} \mbox{$H={\rm diag}\{\hat{H},~1\}$} with
\mbox{$\hat{H}={\rm diag}\{-1,~1\}$,} see also \eqref{Lyap1}.
Namely the corresponding rational functions have {\em at most}
one pole in each open half plane. We present the rational function
along with the corresponding system matrix.
\begin{center}
$~L_{\alpha}=${\mbox{\tiny$\begin{pmatrix}
0&1&1\\ 1&1&0\\ 1&0&0
\end{pmatrix}$}}\quad\quad\quad
$L_{\beta}=${\mbox{\tiny$
\begin{pmatrix}-1&0&0\\~~0&1&1\\~~0&1&1
\end{pmatrix}$}}\quad\quad\quad
$~L_{\gamma}=${\mbox{\tiny$
\begin{pmatrix}-1&0&1\\~~0&1&0\\~~1&0&0\end{pmatrix}$}}

$~\psi_{\alpha}(s)=
\frac{s-1}{s^2-s-1}
\quad\quad\quad\quad\quad
\psi_{\beta}(s)=\frac{s}{s-1}\quad\quad\quad\quad\quad
\psi_{\gamma}(s)=\frac{1}{s+1}$
\vskip 0.3cm

$L_{\delta}=${\mbox{\tiny$\begin{pmatrix}
-1&-1&1\\~~1&~~1&1\\~~1&-1&0
\end{pmatrix}$}}\quad
$~L_{\xi}=${\mbox{\tiny$\begin{pmatrix}
0&0&1\\
0&1&0\\
1&0&1\end{pmatrix}$}}\quad
$~L_{\eta}=${\mbox{\tiny$\begin{pmatrix}-1&-1&2\\-1&~1&1
\\~2&~1&1\end{pmatrix}$}}\quad
$~L_{\theta}=${\mbox{\tiny$\begin{pmatrix}-1&1&-2\\~~1&1&~1\\
-2&1&~1\end{pmatrix}$}}

$\psi_{\delta}(s)=\frac{-4}{s^2}$\quad\quad\quad\quad
$~\psi_{\xi}(s)=\frac{s+1}{s}~$\quad\quad\quad\quad\quad
$~\psi_{\eta}(s)=\psi_{\theta}(s)=\frac{s^2+5s-1}{s^2-2}~.$
\end{center}
Indeed, $~\psi_{\alpha}(s)~$ and
$~\psi_{\eta}(s)=\psi_{\theta}(s)~$ have
a pole in each open half plane.
$~\psi_{\beta}(s)$ has a single pole in $~\C_+$, while
$~\psi_{\gamma}(s)~$ a pole in $~\C_-$. The poles
of
$~\psi_{\delta}(s)$ and $~\psi_{\xi}(s)$
are confined to the imaginary axis.
}
\qed
\end{example}
\vskip 0.2cm

We now point out that unlike to the approach of
\cite[Theorem 2]{FP}, in the PRL framework, see
Proposition \ref{gpLemmaSufficiency}, there is no
restriction on the spectrum of $~L~$ in \eqref{Lyap},
nor on the spectrum of its upper left block, $~A$
in \eqref{Lyap1}. Formally, this follows from Theorem
\ref{Th:Loewy}. Intuitively, restrictions of the form
$~\lambda_j+\lambda_k^*\not=0~$ are
necessary when $~Q$, the right hand side of the Lyapunov
equation is given and a ~{\em unique}~ Hermitian
solution is sought, see e.g. \cite[Corollary 4.4.7]{HJ2}.
In contrast, here we address a genuine Lyapunov
{\em inclusion}. In fact, in each of the cases $~L_{\beta}$,
$~L_{\gamma}$, $~L_{\delta}$, $~L_{\xi}$, $~L_{\eta}~$ and
$~L_{\theta}~$ in Example \ref{InertiaExample},
the $~A~$ matrix (the $n\times n$ upper left block) does
not satisfy this spectral condition.

\section{Convex invertible cones and the
Lyapunov inclusion}
\setcounter{equation}{0}
\label{Sec:6}

In this section we explore some properties of the set
$~\overline{\mathbb L}(H)~$ \eqref{CommonLyap}, to be
used in the sequel in conjunction of the PRL,
see e.g. Theorem \ref{GP(r,nu,p)cic}

We next resort to some background.
Recall that a set of a square matrices is called Convex
Invertible Cone, {\bf cic}~ in short, if it is closed
under positive scaling, summation and inversion. More
precisely, it may include singular elements provided
that the inverse of every non-singular element, belongs
to the set. Thus, the set $~\overline{\mathbb P}~$
is a {\bf cic}. For a more detailed study of
{\bf cic}s see e.g. \cite[Section 2]{CL1},
\cite[Section 2]{CL2} and \cite[Sections 2, 3]{CL4}
Recall that for an Hermitian non-singular matrix $~H~$
we defined in \eqref{CommonLyap} sets of matrices
sharing a common Lyapunov factor, $~{\mathbb L}(H)$,
$~\overline{\mathbb L}(H)$. The following properties
of these sets are fundamental to our discussion.

\begin{theorem}
\label{L(H)}
Let $~H~$ be $~r\times r~$ Hermitian nonsingular, i.e.
\mbox{${\rm inertia}(H)=(\nu,~0,~r-\nu)$}
for some $~\nu\in[0,~r]$.
\begin{itemize}
\item[(i)~~~]{}${\mathbb L}(H)~$ is an invertible cone of
matrices sharing the same inertia as $~H$. It is
a maximal open convex set of nonsingular
matrices \begin{footnote}{For an impressive
particular converse, see
\cite{Ando1}, \cite{Ando2}.}\end{footnote}.

\item[(ii)~~]{}$\overline{\mathbb L}(H)~$ is a closed
invertible cone of matrices with at most $~\nu~$
eigenvalues in $~\C_-~$ and at most $~r-\nu~$ eigenvalues
in $~\C_+$. It is a maximal convex set with this property.

\item[(iii)~]{} $\overline{\mathbb L}(H)^*=\overline{\mathbb L}(H^{-1})$.\\
In particular, $~\overline{\mathbb L}(H)^*=\overline{\mathbb L}(H)$,
if and only if, up to positive
scaling, $~H~$ is an involution $(H^2=I)$.

\item[(iv)~]{} Let $~E~$ be an involution which commutes
with $~H~$ then,\\
$E\overline{\mathbb L}(H)=\overline{\mathbb L}(H)E=
\overline{\mathbb L}(EH)$.

\item[(v)~~]{} For arbitrary Hermitian involution $~E$,
$~E\overline{\mathbb L}(E)=\overline{\mathbb L}(E)E=
\overline{\mathbb L}(I)$.
%
\end{itemize}
\end{theorem}

{\bf Proof}
(i) See \cite[Proposition 3.7]{CL1}, \cite[Observation 2.3.3]{CL4}.

(ii) Part of the claim appeared in \cite[Proposition 2.4(i)]{CL2}.
To establish the inertia property assume first that $~H=I_r$. By
Theorem \ref{Th:Loewy} all matrices in $~\overline{\mathbb L}(I)~$
have inertia $~(0,~r-\pi,~\pi)~$ for some $~\pi\in[0,~r]$. Let
$~B~$ be an $~r\times r~$ matrix not in
$~\overline{\mathbb L}(I)$, namely the Hermitian part of $~B~$
has a negative eigenvalue, i.e. $~\min\limits_{k=1,~\ldots~,~r}
\lambda_k(B+B^*)=-\beta~$ for some $~\beta>0$. Take now
$~A=\frac{\beta}{4}I+\frac{1}{2}(B^*-B)$. It is straightforward
to verify that $~A\in{\mathbb L}(I)$, but
$~A+B=\frac{\beta}{4}I+\frac{1}{2}(B+B^*)$. Thus,
$~\min\limits_{k=1,~\ldots~,~r}\lambda_k(A+B)=-\frac{\beta}{4}$,
i.e. in contrast to all elements of $~\overline{\mathbb L}(I)$,
the matrix $~A+B~$ has eigenvalues in $~\C_-$ so this part of the
claim is established for $~H=I$. To address a general Hermitian
non-singular matrix $~H$, exploit the fact, see Observation
\ref{congruence} below, that one can always find a non-singular
$~V~$ and an involution $~E_{\nu}~$ \eqref{E} so that
$~E_{\nu}V^{-1}\overline{\mathbb L}(H)V =V^{-1}\overline{\mathbb
L}(H)VE_{\nu}=\overline{\mathbb L}(I)$, so the claim is
established.

(iii) The claim for arbitrary $~H~$ appeared in
\cite[Equation (3.6)]{CL1}. If $~H~$ is a scaled involution, i.e.
$~H^{-1}=\alpha{H}~$ for some $~\alpha>0$, the claim is
follows from multiplying $~\overline{\mathbb L}(H)~$ in
\eqref{CommonLyap} by $~H^{-1}~$ from both sides.

For the other direction assume that $~\left(\overline{\mathbb
L}(H)\right)^*= \overline{\mathbb L}(H)$. This relation is
invariant under unitary similarity. Thus, without loss of
generality assume that $~H~$ is (real non-zero) diagonal, say
$~H={\rm diag}\{h_1~,~\ldots~,~h_r\}$. Take now $~A~$ to be equal
to $~H$, except a single non-zero off diagonal element $~x~$ at
the location $~jk$, where $~j>k$, i.e. $~A~$ is lower triangular.
A straightforward calculation shows that $~A\in\overline{\mathbb
L}(H)~$ is equivalent to $~2|h_j|\geq |x|$. By assumption, also
$~A^*\in\overline{\mathbb L}(H)$, which in turn is equivalent to
$~2|h_k|\geq |x|$. Thus, $~h_j=\pm{h_k}~$ and since true for
all $j,k$ this claim is established.

(iv) This claim is proved for $~{\mathbb L}(H)$, in
\cite[Lemma 3.6]{CL1}, see also \cite[Proposition 3.2.2]{CL4}.
The case $~\overline{\mathbb L}(H)~$ is similar and thus omitted.
Item (v) is in fact a special case of item (iv).
%
%
%
\qed
\vskip 0.2cm

Note that item (iv) in Theorem \ref{L(H)} says that if
$~L\in\overline{\mathbb L}(H)$, for some Hermitian non-singular
$~H$, then both matrices $~EL~$ and $~LE~$ belong to
$~\overline{\mathbb L}(EH)$, whenever $~E~$ is an involution
which commutes with $~H$. Using this, along with
item (v) in Theorem \ref{L(H)}, we state the following.

\begin{corollary}\label{H(E)}
For natural $~n, p~$ and $~\nu\in[1, n]~$ let us denote
$~{H}_o={\rm diag}\{I_{\nu},~-I_{n-\nu},~I_p\}$,
$~H_1={\rm diag}\{-I_{\nu},~I_{n-\nu},~I_p\}~$ and
$~H_2={\rm diag}\{-I_n,~I_p\}$.
Then,
\[
\begin{matrix}
H_o\overline{\mathbb L}(H_1)&=&
\overline{\mathbb L}(H_1)H_o&=&
\overline{\mathbb L}(H_2),\\
H_o\overline{\mathbb L}(H_2)&=&
\overline{\mathbb L}(H_2)H_o&=&
\overline{\mathbb L}(H_1).
\end{matrix}
\]
In addition, $~H_j\overline{\mathbb L}(H_j)=
\overline{\mathbb L}(H_j)H_j=
\overline{\mathbb L}(I_{n+p})\quad\quad j=1, 2.$

Namely, there is one-to-one correspondence between the sets
$~\overline{\mathbb L}(H_1)$,
$~\overline{\mathbb L}(H_2)~$ and $~\overline{\mathbb L}(I_{n+p})$.
\end{corollary}
\vskip 0.2cm

We shall find it convenient to introduce the following
notation for $~l\times l~$ signature matrices,
\begin{equation}\label{E}
\begin{matrix}
E_{\nu, l}:={\rm diag}\{-I_{\nu}~,~I_{l-\nu}\}
&~&~&\nu\in[0,~l].\end{matrix}
\end{equation}
Whenever the dimension $~l~$ is evident from the context
we shall simply write $~E_{\nu}$.
\vskip 0.2cm

We can now cite the following known facts
see e.g. \cite[Lemma 3.4]{CL1}.

\begin{observation}\label{congruence}
Consider the $~r\times r~$ nonsingular matrices $~V~$
and $~H=H^*$. Then the following relations
hold in \eqref{CommonLyap},
$$\begin{matrix}
(a)~V^{-1}{\mathbb L}(H)V={\mathbb L}(V^*HV)
&~&~&
(b)~V^{-1}\overline{\mathbb L}(H)V=
\overline{\mathbb L}(V^*HV).\end{matrix}$$
In particular, if $~{\rm inertia}(H)=(\nu,~0,~n-\nu)$,
for some $~\nu\in[0,~r]$, $~V~$ may be chosen so that,
$$\begin{matrix}
(a)~V^{-1}{\mathbb L}(H)V={\mathbb L}(E_{\nu})
&~&
(b)~V^{-1}\overline{\mathbb L}(H)V=
\overline{\mathbb L}(E_{\nu}),
\end{matrix}$$
where $~E_{\nu}~$ is as in \eqref{E}.
%
\end{observation}
\vskip 0.2cm

From Observation \ref{congruence} it follows that, up to
similarity,
$~\bigcup\limits_{\nu=0}^r\overline{\mathbb L}(E_{\nu, r})~$
covers all $~r\times r~$ matrices. We next point out that
technically this is not a proper partitioning. A
straightforward substitution in \eqref{CommonLyap} with both
$~H=E_{\nu}~$ and $~E_{\nu+\eta}$, reveals that these sets
are ``nearly" distinct.

\begin{corollary}\label{L(E2)}
Let $~\nu\geq 0~$ and $~\eta\geq 1~$ be so that $~r\geq\nu+\eta$
and $~E_{\nu}~$ is as in \eqref{E}. If $~L\in\{\overline{\mathbb
L} (E_{\nu})\bigcap\overline{\mathbb L}(E_{\nu+ \eta})\}~$ then
$~L=${\mbox{\tiny$\begin{pmatrix} -Q_1+T_1&0&K-R\\ 0&T_3&0\\
K^*&0&Q_2+T_2\end{pmatrix}$}} where $~T_1, T_2,~T_3~$ are
skew-Hermitian matrices of dimensions $~\nu\times\nu$,
$~\eta\times\eta~$ and $~(r-\nu-\eta)\times(r-\nu-\eta)$,
respectively, $~K\in\C^{\nu\times(r-\nu-\eta)}$ is arbitrary
and the other blocks are so that the
$~(r-\eta)\times(r-\eta)~$ matrix
{\mbox{\tiny$\begin{pmatrix}2Q_1&R\\ R^*&2Q_2
\end{pmatrix}$}} is positive semi-definite.
\end{corollary}
\vskip 0.2cm

Observation \ref{congruence} also suggests that for every
non-singular Hermitian $~H~$ the set $~\overline{\mathbb L}(H)~$
is isomorphic to $~\overline{\mathbb L}(I)$. We conclude this
section by showing that in turn, one can describe
$~\overline{\mathbb L}(I)~$ as a sum of two {\bf cic}s:
$~\overline{\mathbb P}~$ and $~{\mathbb T}$, the set of
skew-Hermitian matrices, see e.g. \cite[Proposition
3.2.5(ii)]{CL4}. Furthermore, each may be described by the convex
hull of its extreme directions. To this end, we need to introduce
$~\mathbb{OP}$, the set of rank one orthogonal projections, i.e.
\[
\mathbb{OP}=\{ \pi=xx^*~:~
\begin{matrix}x\in\C^r&~&x^*x=1\end{matrix}~\}.
\]
\begin{observation}\label{L(I)}

\noindent
I. Let $~E~$ be as in \eqref{E}, then
\[
\begin{matrix}
E\overline{\mathbb L}(I)&=&\overline{\mathbb L}(I)E&=&
\overline{\mathbb L}(E),\\
E{\mathbb L}(I)&=&{\mathbb L}(I)E&=&{\mathbb L}(E).
\end{matrix}
\]
II. \cite[Proposition 3.2.5]{CL4}
$$
{\mathbb L}(I)={\mathbb P}+{\mathbb T}\quad\quad\quad\quad
\overline{\mathbb L}(I)=\overline{\mathbb P}+{\mathbb T}.
$$
III. The sets $~\overline{\mathbb P}$ and
$~{\mathbb T}~$ may be constructed from orthogonal
projection. Indeed,
$$\overline{\mathbb P}=\left\{ \sum\limits_{j=1}^r\alpha_j\pi_j~:~
\begin{matrix}\alpha_j\geq 0&~&\pi_j\in{\mathbb O}{\mathbb P}
\end{matrix}\right\},$$
where $~\pi_1,~\ldots~,~\pi_r~$ are all distinct.
Similarly,
$${\mathbb T}=\left\{ i\sum\limits_{j=1}^r\rho_j\pi_j~:~
\begin{matrix}&
\rho_j\in\R&~&\pi_j\in{\mathbb O}{\mathbb P}\end{matrix}
\right\}$$
where $~\pi_1,~\ldots~,~\pi_r~$ are all distinct.
\end{observation}

To summarize, for arbitrary $~E~$ there is a one-to-one
correspondence between the sets $~\overline{\mathbb L}(E)~$
and $~\overline{\mathbb L}(I)$. Thus, it suffices to
construct the latter set. Indeed,
$~\overline{\mathbb L}(I)=\overline{\mathbb P}+
{\mathbb T}$. Now,
$~P\in\overline{\mathbb P}~$ can always be parameterized
by non-negative scalars $~\alpha_1,~\ldots~,~\alpha_r~$
and $~r~$ distinct points on the $~\|~\|_2~$ unit sphere.
Similarly, $~T\in\mathbb{T}~$ can always be parameterized
by real scalars $~\rho_1,~\ldots~,~\rho_r~$
and $~r~$ distinct points on the $~\|~\|_2~$ unit sphere.
\vskip 0.2cm

In fact, parameterization of a point on $~\|~\|_2~$ unit
sphere can be further simplified.  Note that
$~\pi\in\mathbb{OP}~$ may be identified with a point on the
$~\|~\|_2~$ unit sphere \mbox{$\{~x\in\C^r~:~x^*x=1~\}$}.
Now, through polar coordinates there is
a one-to-one correspondence between this unit sphere and
\mbox{$\{~y\in\R^{r(r-1)}~:~2\pi>\| y\|_{\infty}\}$.}
For example for $~r=3$,
$~x=${\mbox{\tiny$\begin{pmatrix}
\cos(\theta_1)\cos(\theta_2)e^{i\eta_1}\\
\sin(\theta_1)\cos(\theta_2)e^{i\eta_2}\\
\sin(\theta_2)e^{i\eta_3}
\end{pmatrix}$}} $\theta_j,\eta_k\in[0,~2\pi)$, $~j=1, 2$ and
$~k=1,2,3$.
Thus elements in $~\mathbb{OP}~$ can be parameterized by
points in the real ``box" $~[0,~2\pi)^{r(r-1)}$.
\vskip 0.2cm

In the next section we identify the matrix {\bf cic}
$~\overline{\mathbb L}(H)~$ with the set of system
matrices associated with a subset of $~\mathcal {GP}~$
functions, see Theorem \ref{GP(r,nu,p)cic} below.
These (not necessarily minimal) realizations, cover
all $~\mathcal {GP}~$ functions with no pole at infinity.

\section{convex invertible cones of realizations of generalized
positive functions of prescribed parameters}
\setcounter{equation}{0}
\label{Sec:7}

We start with a couple of known related results. First,
recall that in \eqref{DefGP} we have already described
$~\mathcal{GP}~$ functions as map from $~i\R~$ to
$~\overline{\mathbb L}(I)$. We now formalize this fact.

\begin{observation}\label{GP-L(I)}
One can view $p\times p$-valued $~\mathcal{GP}~$ functions as
a {\bf cic} of rational functions (almost everywhere)
analytically mapping $~i\R~$ to $~\overline{\mathbb L}(I_p)$,
\eqref{CommonLyap}, and $~{\mathcal P}~$ as a sub{\bf cic}
mapping $~\C_+~$ to $~\overline{\mathbb L}(I_p)$.
\end{observation}
\vskip 0.2cm

In \cite[Section 4.4]{CL4} analogies
were drawn between the set $~\mathcal{P}~$ of scalar
rational functions and the matrix set
$~{\mathbb L}(H)~$ for $~H\in{\mathbb P}$. We now
use Observation \ref{GP-L(I)} to introduce an
analogy between the sets $~\mathcal{GP}~$ and
$~{\mathbb L}(I)$. Recall that
$$
V\overline{\mathbb L}(I)V^*\subset\overline{ \mathbb L}(I)
\quad\quad\quad V\in\C^{n\times n},
$$
(and if in addition $~V~$ is
non-singular, then, $~V{\mathbb L}(I)V^*\subset{\mathbb
L}(I)$). Similarly,
$$
F{\mathcal G}{\mathcal P}F^{\#}\subset\mathcal{GP}
\quad\quad\quad F\in{\mathcal F}.
$$
Recall also that
\vskip 0.2cm

As another association between the matrix {\bf cic}
$~\overline{\mathbb L}(H)~$ and a  {\bf cic} of
rational functions, we cite the following. Consider the
set of all rational functions of the form $~C(sI-A)^{-1}B$
admitting~ {\em balanced realization}, i.e.
$$HA^*+AH=BB^*\quad\quad\quad HA+A^*H=C^*C,$$
for some non-singular Hermitian $~H$. Each of these sets
forms a {\bf cic} of state space realizations. In particular,
this allows for simultaneous model order reduction of
uncertain systems, see \cite[Section 5]{CL2} for details.
\vskip 0.2cm

We now identify the matrix {\bf cic} $~\overline{\mathbb L}(H)$,
\eqref{CommonLyap}, with a {\bf cic} of system matrices,
\eqref{StateSpace}, associated with a subset of
$~\mathcal {GP}~$ functions. As a motivation recall that
the set of $~p\times p~$ $~\mathcal{GP}~$ functions is a
{\bf cic} of rational functions. However, the McMillan
degree of a sum, is roughly the sum of the McMillan
degree of the original functions. Now, if one considers a
pair of $~p\times p~$ $~\mathcal{GP}~$ functions, of McMillan
degree at most $~n$, admitting state space realizations of the
form \eqref{StateSpace}, the sum of the respective system
matrices is associated with a $~p\times p~$ rational
function of McMillan degree at most $~n$. However, this
``sum" function may be not generalized positive. This is
illustrated next.

\begin{example}\label{ConvRealiz}
{\rm
For simplicity, consider two (minimal) state space
realizations of the scalar function
$~\psi_{\epsilon}(s)=\psi_{\delta}(s)=\frac{-4}{s^2}~$
from Example \ref{InertiaExample}:
$L_{\delta}=${\mbox{\tiny$\begin{pmatrix}
-1&-1&1\\~~1&~~1&1\\~~1&-1&0\end{pmatrix}$}}
and
$~L_{\epsilon}=${\mbox{\tiny$\begin{pmatrix}~~1&-1&~~1
\\~~1&-1&-1\\-1&-1&~~0\end{pmatrix}$}}
(note that \mbox{$L_{\epsilon}=V^{-1}L_{\delta}V$} with
$~V={\rm diag}\{\hat{V}~,~~1\}~$ where $~\hat{V}=$
{\mbox{\tiny$\begin{pmatrix}0&-1\\ 1&~~0\end{pmatrix}$}}
).
Let now $~L_{\zeta}~$ be a convex combination of these two
realizations, i.e.
\mbox{$L_{\zeta}=\frac{1}{2}(L_{\delta}+L_{\epsilon})$}
$=${\mbox{\tiny$\begin{pmatrix}0&-1&1\\1&0&0\\0&-1&0
\end{pmatrix}$}}. This is a (minimal) state space
realization of $~\psi_{\zeta}(s)=\frac{-1}
{s^2+1}_{|_{s=0}}=-1$, which is not generalized
positive.}\qed
\end{example}
\vskip 0.2cm

{\bf Coordinates transformation}

Following Proposition \ref{gpLemmaSufficiency}, we have
focused our attention on a special case of the set
$~\overline{\mathbb L}(H)$, where \mbox{$r=n+p$} and the
Lyapunov factor $~H~$ is block diagonal of the form
$~H=\rm{diag}\{\hat{H},~I_p\}$, where $~\hat{H}~$
is $~n\times n~$ Hermitian non-singular.

From Observation \ref{congruence} it follows that up to
similarity over
$~\overline{\mathbb L}(H)_{|_{H=\rm{diag}\{\hat{H},~I_p\}}}$,
one can confine the discussion to the case where in
addition in \eqref{CommonLyap}
$~\hat{H}=E_{\nu, n}$, see \eqref{E}.
\vskip 0.2cm

Recall that coordinates transformation means that whenever
$~V={\rm diag}\{\hat{V}, I_p\}~$ with $~\hat{V}~$ $~n\times n~$
nonsingular, with $~L~$ in \eqref{StateSpace}, $~V^{-1}LV~$ is
another state space realizations of the same rational function
$~\Psi(s)$. Thus, taking in \eqref{E} $~l=r=n+p$, without
loss of generality we can focus on sets of $~(n+p)\times(n+p)~$
matrices of the form
\[
\overline{\mathbb L}(H)\quad\quad\quad
H={\rm diag}\{E_{\nu, n},~I_p\}\quad\quad\quad\nu\in[0, n],
\]
partitioned as
{\mbox{\tiny$\begin{pmatrix}A&B\\ C&D\end{pmatrix}$}}. We
now find it convenient to denote by
\begin{equation}\label{voltaire}
\mathcal{GP}(r,~\nu,~p)\quad\quad
r\geq2\quad p\in[1,~r-1] \quad \nu\in[0,~r-p]
\end{equation}
the set of all $~p\times p$-valued rational functions obtained
by \eqref{StateSpace} and \eqref{Lyap} (recall, $~A~$ is
$~(r-p)\times (r-p)~$). From item (ii) of
Theorem \ref{L(H)} we have the following:

\begin{theorem}\label{GP(r,nu,p)cic}
Given a set $~\mathcal{GP}(r,~\nu,~p)~$ as described in
\eqref{voltaire}. The set of all the corresponding
state space matrices forms a closed invertible cone of
realizations of functions with at most $~\nu~$ poles in
$~\C_-~$ and at most $~r-p-\nu~$ poles in $~\C_+$. Moreover,
this is a maximal convex realization set with this property.
\end{theorem}

Maximality is in the sense described in proof of
item (ii) of Theorem \ref{L(H)}: Recall that each of the
realization matrices $~L~$ associated with functions in
$~\mathcal{GP}(r,~\nu,~p)~$ has at most $~\nu~$ eigenvalues in
$~\C_-~$ and at most $~r-\nu~$ eigenvalues in $~\C_+$. Now
if $~\tilde{L}~$ is a $~r\times r~$ realization matrix
associated of a rational function
$~\tilde{\psi}\not\in\mathcal{GP}(r,~\nu,~p)$,
then one can always find $~L~$ so that $~L+\tilde{L}~$ has
$~\nu+1~$ in $~\C_-~$ or $~r-\nu+1~$ eigenvalues in $~\C_+$.
\vskip 0.2cm

Before proceeding, we now find it convenient to denote by
$~\mathcal{GP}_{\rm min}(r,~\nu,~p)~$ the subset of
$~\mathcal{GP}(r,~\nu,~p)~$ where in addition the realization
is minimal (i.e. the McMillan degree $~q~$ is equal to $~r-p$).
Under this terminology the necessity part of the PRL,
i.e. Lemma \ref{gpLemma}, and Proposition
\ref{gpLemmaNecessity1}, is restricted to elements in
$~\mathcal{GP}_{\rm min}(r,~\nu,~p)$.
\vskip 0.2cm

\begin{example}\label{GP(r,nu,p)}
{\rm
We next illustrate some aspects of the correspondence
between the matrix set $~\overline{\mathbb L}(E_{\nu})~$
and the set of rational functions $~\mathcal{GP}(r,~\nu,~p)$,
along with its subset $~\mathcal{GP}_{\rm min}(r,~\nu,~p)$.
As before, we concentrate on the case where $~r=3$,
$~\nu=1~$ and $~p=1$.

\begin{itemize}
\item[(i)~~~~]{}To illustrate Theorem \ref{GP(r,nu,p)cic}
the functions
$~\psi_{\alpha}(s)$,
$~\psi_{\beta}(s)$, $~\psi_{\gamma}(s)$, $~\psi_{\delta}(s),$
$~\psi_{\xi}(s)~$ and $~\psi_{\eta}(s)~$ from Example
\ref{InertiaExample} are all in $~\mathcal{GP}(3, 1, 1)$.
\vskip 0.2cm

\item[(ii)~~~]{}Starting from $~L$, one may obtain, minimal and
non-minimal realizations. $~\psi_{\beta}(s)$,
$~\psi_{\gamma}(s)~$ and $~\psi_{\xi}(s)$, from Example
\ref{InertiaExample} are in
\mbox{$\mathcal{GP}(3, 1, 1)\smallsetminus\mathcal{GP}_{\rm min}(3, 1, 1)$.}
In contrast, starting from rational functions,
$~\psi_{\gamma}(s)~$ and $~\psi_{\xi}(s)~$
can
be written as part of $~\mathcal{GP}_{\rm min}(2,~1,~1)$,
i.e. $~\hat{L}_{\gamma}=${\mbox{\tiny$\begin{pmatrix}-1&0\\~~1&0
\end{pmatrix}$}} and
\mbox{$\hat{L}_{\gamma}^{-1}=\hat{L}_{\xi}=$
{\mbox{\tiny$\begin{pmatrix}0&1\\ 1&1\end{pmatrix}$}}} are the
corresponding system matrices. This now conforms with
Proposition \ref{PosEquiv1} part (I) below noting that they
are positive.

\item[(iii)~~]{}${\mathbb L}(E_1)~$ is a convex invertible
cone, {\bf cic}. Samples of it are all realization presented
in Example \ref{InertiaExample}. In particular,
\mbox{$L_{\xi}=L_{\gamma}^{-1}$} and
\mbox{$L_{\beta}=\frac{1}{2}(L_{\eta}+L_{\theta})$.}\\
In contrast, in Example \ref{ConvRealiz}, $~L_{\epsilon}~$ is
constructed from $~L_{\delta}~$ so that they do not share a
common Lyapunov factor, of the form
$~{\rm diag}\{\hat{H}~,~1\}$, to the LMI in \eqref{Lyap}.

\item[(iv)~~]{}Even in $~\mathcal{GP}_{\rm min}(3, 1, 1)$
elements may not have a pole in each open half plane, see
e.g. $~\psi_{\delta}(s)~$ in Example \ref{InertiaExample}.

\item[(v)~~~]{}In contrast to the set
$~\mathcal{GP}(r,~\nu,~p)$, the set of realization
matrices associated with the subset
$~\mathcal{GP}_{\rm min}(r,~\nu,~p)~$, is not convex.
Recall that on the one hand,
$~\psi_{\eta}(s)=\psi_{\theta}(s)$ is in
$~\mathcal{GP}_{\rm min}(3, 1, 1)$, i.e. $~L_{\eta}$,
$L_{\theta}~$ are the respective minimal realizations.
However, $~L_{\beta}=\frac{1}{2}(L_{\eta}+L_{\theta})~$
is a non-minimal realization of $~\psi_{\beta}(s)$.
\end{itemize}
\qed
}
\end{example}
\vskip 0.2cm

We now use the Lyapunov inclusion formulation to establish
a strong link between positive and
generalized positive rational functions.

\begin{lemma}\label{LemmaPosEquiv}
Let $~r,~\nu,~p~$ with $~r\geq 2$, $p\in[1,~r-1]$ and
$~\nu\in [0,~r-p]$ be given and let
$~J={\rm diag}\{I_{\nu},~-I_{r-p-\nu}\}$
(i.e. $J=-E_{\nu, r-p}$).

Let {\mbox{\tiny$\begin{pmatrix}A&B\\ C&D\end{pmatrix}$}}
be a system matrix associated a function in
$~\mathcal{GP}(r,~\nu,~p)~$ and let
{\mbox{\tiny$\begin{pmatrix}\hat{A}&\hat{B}\\ \hat{C}&\hat{D}
\end{pmatrix}$}} be associated with a $~p\times p$-valued
positive rational function whose realization is of dimension $~r-p$.

Then, the following is true.
\begin{itemize}
\item[(I)~]{}
{\mbox{\tiny$\begin{pmatrix}JA&JB\\ C&D\end{pmatrix}$}}
and
{\mbox{\tiny$\begin{pmatrix}AJ&B\\ CJ&D\end{pmatrix}$}}
are two $~r-p~$ dimensional realizations of the same
$~p\times p$-valued positive rational function.
\vskip 0.1cm

\item[(I)~]{}
{\mbox{\tiny$\begin{pmatrix}J\hat{A}&J\hat{B}\\ \hat{C}&\hat{D}
\end{pmatrix}$}} and
{\mbox{\tiny$\begin{pmatrix}\hat{A}J&\hat{B}\\ \hat{C}J&\hat{D}
\end{pmatrix}$}} are two realizations of the same function in
$~\mathcal{GP}(r,~\nu,~p)$.
\end{itemize}
\end{lemma}

{\bf Proof :}\quad
To see that these are two realization of the same rational
function, recall that for arbitrary involution $~J$, one
has that,
\begin{center}
$C(sI_n-JA)^{-1}JB+D=
C(sJ-A)^{-1}B+D=
CJ(sI_n-AJ)^{-1}B+D.$
\end{center}
For the equivalence, the claim follows from Proposition
\ref{gpLemmaSufficiency} and Corollary \ref{H(E)} with
$~n=r-p~$ and $~H_o=J$.
\qed
\vskip 0.2cm

This can be formalized in the following stating that
for given $~r>p\geq 1$, the set of $~p\times p$-valued
rational generalized positive functions whose realization
is of dimension $~r-p$, can be described as
$~r+1-p~$ replicas of its subset of positive functions.

\begin{proposition}\label{PosEquiv1}
Let $~r, \nu, p~$ be arbitrary.
\begin{itemize}
\item[(I)~~]{} If $~r=\nu+p~$ and $~\nu\geq 1$ then,
\[
\mathcal{GP}(r,~\nu,~p)_{|_{r=\nu+p}}\subset
\mathcal{P}.
\]
\item[(II)~]{} The set $~\mathcal{GP}(r,~\nu,~p)~$
is state-space equivalent to the set of
positive functions $~\mathcal{GP}(r,~r-p,~p)$.

\item[(III)]{}For given $~r>p\geq 1$, the set of
$~p\times p$-valued rational generalized positive
functions whose realization is of
dimension $~r-p$, can be described as
\[
\bigcup\limits_{\nu=0}^{r-p}\mathcal{GP}(r,~\nu,~p).
\]
\end{itemize}
\end{proposition}

To see that the inclusion in item (I) of Proposition \ref{PosEquiv1}
is strict, take for example $~\psi_{\gamma}(s)~$ from Example
\ref{InertiaExample} where $~r=3$, $~\nu=1~$ and $~p=1$, i.e. it is a
positive function in \mbox{$\mathcal{GP}(3,~1,~1)\smallsetminus
\mathcal{GP}_{\rm min}(3,~1,~1)$,} but neither belongs to
$~\mathcal{GP}(3,~1,~2)$ nor to $~\mathcal{GP}(3,~2,~1)$. See also
item (ii) in Example \ref{GP(r,nu,p)}.
\vskip 0.2cm

Strictly speaking the sets $~\mathcal{GP}(r,~\nu,~p)~$ do not
offer a partitioning of generalized positive function to distinct
subsets. Namely, from Corollary \ref{L(E2)} is follows that
$~\mathcal{GP}(r,~\nu_1,~p)\cap\mathcal{GP}(r,~\nu_2,~p)\not=\emptyset$.
For example, the system matrix
$~L=${\mbox{\tiny$\begin{pmatrix}0&~0&0\\ 0&~a&b\\
0&-b&0\end{pmatrix}$}} with $a\geq 0$,
is a realization of $\psi=\frac{b^2}{a-s}$. As
\mbox{$L\in\mathbb{L}(I)\cap\mathbb{L}(E_1)$,} this
$~\psi~$ belongs to
$~\mathcal{GP}(3, 0, 1)\cap\mathcal{GP}(3, 1, 1)$.
\vskip 0.2cm

Observation \ref{L(I)} offered an easy-to-compue construction of
$~\overline{\mathbb L}(I)~$ and subsequently of
$~\overline{\mathbb L}(E)$, where $~E=E_{\nu, r}~$
and $~\nu\in[0, r-p]~$ is arbitrary. Theorem
\ref{GP(r,nu,p)cic} and Eq. \eqref{voltaire} identify the matrix set
$~\overline{\mathbb L}(E)~$ with $~\mathcal{GP}(r,~\nu,~p)$. Thus,
through item (III) of Proposition \ref{PosEquiv1} we here offer
a construction of the set of all $~p\times p$-valued rational
$~\mathcal{GP}~$ functions whose realization is of dimension $~r-p$,
where $~r>p\geq 1$, are arbitrary.
\vskip 0.2cm

LMI control theory and Matlab based LMI procedures, were
developed essentially for $~\mathcal{P}~$ functions, see e.g.
\cite{BGFB}. However, as stated in Proposition \ref{PosEquiv1},
from realization point of view,
for arbitrary $r\geq 2$, $p\in[1,~r-1]$ and $\nu\in[0,~p-r]$,
the set $~\mathcal{GP}(r,~\nu,~p)~$ functions is equivalent
to the set $~\mathcal{GP}(r,~r-p,~p)~$ positive functions.
It is thus of interest to try to adapt some of the LMI theory
to $~\mathcal{GP}~$ functions.
\vskip 0.2cm

To further motivate the introduction of the set
$~\mathcal{GP}(r, \nu, p)$, in the next section we
provide a
classical control interpretation of it.

\section{turning a function generalized positive
through static output feedback}
\setcounter{equation}{0}
\label{Sec:8}

Let $~F(s)~$ be a $~p\times p$-valued rational function
vanishing at infinity $(\lim\limits_{s\rightarrow\infty}F(s)=0)$.
Thus, it admits a space realization of the form
$~F(s)=C(sI-A)^{-1}B$,
\begin{equation}\label{OpenLoop}
\dot{x}=Ax+Bu,\quad\quad\quad y=Cx,
\end{equation}
where $~u, y~$ are $p-$dimensional and
$~x~$ is $~n-$dimensional\begin{footnote}
{The case $D\not=0$ is more involved and thus omitted.}
\end{footnote}.
Using \eqref{StateSpace}, the corresponding
system matrix is $~L=${\mbox{\tiny$\begin{pmatrix}A&B\\
C&0_p\end{pmatrix}$}}.
Assume that a static output feedback is applied i.e.
\begin{equation}\label{F}
u=Ky+u',
\end{equation}
where $~K~$ is a constant $~p\times p~$ matrix and
$~u'~$ is an auxiliary input. Thus, the resulting
closed loop system (mapping $~u'~$ to $~y$) is
\begin{equation}\label{Fcl}
F_{\rm cl}(s):=C(sI-A_{\rm cl})^{-1}B,\quad\quad
A_{\rm cl}:=A+BKC.
\end{equation}
The associated system matrix is, $~L_{\rm
cl}=${\mbox{\tiny$\begin{pmatrix}A_{\rm cl} &B\\
C&0_p\end{pmatrix}$}}. Robust stabilization of a system of the
form \eqref{OpenLoop} by static output feedback \eqref{F} is
quite known, see e.g. \cite[Theorem 3.2]{deSX}, \cite[Section
3]{ElGOAiR} \cite{HL,HSKu} and \cite{PFA}. It is well known that
positive functions are closely related with robust stability, see
e.g. \cite{NT}, \cite{F}. We here introduce a characterization of
all systems which may be turned generalized positive through
static output feedback.

\begin{proposition}\label{OutputFeedback}
Let $~F(s)~$ be a $~p\times p$-valued rational function
vanishing at infinity $(\lim\limits_{s\rightarrow\infty}F(s)=0)$
with a state space realization as in \eqref{OpenLoop}.

\begin{itemize}
\item[(I)~~~~]{}There exists a static output feedback,
of the form of \eqref{F}, rendering the closed loop
system, $~F_{\rm cl}(s)~$ in \eqref{Fcl} generalized
positive, if and only if, there exists a nonsingular
Hermitian $~\hat{H}\in\C^{n\times n}$ so that
the open loop system \eqref{OpenLoop} satisfies,

(a)~~$C=-B^*\hat{H}$,

(b)~~$v^*(A\hat{H}^{-1}+\hat{H}^{-1}A^*)v\geq 0~$ for
all vector $~v~$ in the null-space of $~B^*$.

\item[(II)~~~]{} Let $~r, \nu, p$, where $~\nu\in[0,~r-p]~$
and $~p\in[1,~r-\nu]$ be given. There exists a static output
feedback, of the form of \eqref{F}, so that $~F_{\rm cl}(s)~$
in \eqref{Fcl} belongs $~\mathcal{GP}(r, \nu, p)$, if and
only if, up to coordinate transformation,
the open loop system \eqref{OpenLoop} satisfies,

(a)~~$C=-B^*E_{\nu}$,

(b)~~$v^*(AE_{\nu}+E_{\nu}A^*)v\geq 0~$ for
all vector $~v~$ in the null-space of $~B^*$.

\item[(III)~~]{} Let $~r, \nu, p$, where $~\nu\in[0,~r-p]~$ and
$~p\in[1,~r-\nu]$ be given. The set $~\mathcal{GP}(r, \nu, p)~$
is invariant under static output feedback \eqref{F}
whenever $~K\in\overline{\mathbb L}(-I_p)$.

\item[(IV)~~]{} If in part (I) $~-\hat{H}\in{\mathbb P}_n~$
or in parts (II), (III), $~\nu=r-p$, then
the resulting closed loop
function $~F_{\rm cl}$ is in $\mathcal{P}$.
\end{itemize}
\end{proposition}

{\bf Proof~:}\quad
(I)~~ Following Proposition \ref{gpLemmaSufficiency}
denote $~H={\rm diag}\{\hat{H},~I_p\}~$ with
\mbox{$\hat{H}\in\C^{n\times n}$}
Hermitian nonsingular,
$~W:=HL_{\rm cl}+L_{\rm cl}^*H~$ and recall
that having $~F_{\rm cl}\in\mathcal{GP}~$ is
equivalent to $~W\in\overline{\mathbb P}_r$
\begin{center}
$W=${\mbox{\tiny$\begin{pmatrix}
\hat{H}A+A^*\hat{H}&~C^*+\hat{H}B\\ C+B^*\hat{H}&0_p
\end{pmatrix}$}}$+$
{\mbox{\tiny$\begin{pmatrix}\hat{H}BKC+(BKC)^*\hat{H}&0\\
0&0_p\end{pmatrix}$}}.
\end{center}
Thus, having $~W\in\overline{\mathbb P}_r$, implies
condition (a) in the claim. Substituting back
one has that
\begin{center}
$W=${\mbox{\tiny$\begin{pmatrix}
\hat{H}A+A^*\hat{H}&0\\ 0&0_p\end{pmatrix}$}}$+$
{\mbox{\tiny$\begin{pmatrix}-
\hat{H}B(K+K^*)B^*\hat{H}&0\\ 0&0_p\end{pmatrix}$}}
$=${\mbox{\tiny$\begin{pmatrix}\hat{H}\hat{W}
\hat{H}&0\\ 0&0_p\end{pmatrix}$}}
\end{center}
where $~\hat{W}=A\hat{H}^{-1}+\hat{H}^{-1}A^*-B(K+K^*)B^*$.
Now, $~W\in\overline{\mathbb P}_r~$ is equivalent to
$~\hat{W}\in\overline{\mathbb P}_n$ (recall $r=n+p)$.
Next, whenever $~v^*\hat{B}\not=0$, one can take
$~K\in\overline{\mathbb L}(-I_p)~$ ``sufficiently
large", so that $~v^*\hat{W}v\geq 0$. Now, if
$~v^*\hat{B}=0$, condition (b) guarantees that
$~v^*\hat{W}v\geq 0$, so this part of the
claim is established.

(II)~~ One can always write $~\hat{H}=V^*E_{\nu}V~$ for some
non-singular $~V$. The statement follows from applying the
corresponding transformation of coordinates to the
realization of $~F(s)$.

(III)~~ This follows from the proof of part (I) together
with Proposition \ref{gpLemmaSufficiency}, noting that by
assumption $~\hat{H}A+A^*\hat{H}\in\overline{\mathbb P}_n$,
which is equivalent to
$~A\hat{H}^{-1}+\hat{H}^{-1}A^*\in\overline{\mathbb P}_n$.

(IV)~ See part (I) of Proposition \ref{PosEquiv1},
so the claim is established.
\qed
\vskip 0.2cm
%
\vskip 0.2cm

\begin{example}\label{Feedback}
{\rm

We here illustrate Proposition \ref{OutputFeedback}
by examining sets of rational functions
$~\mathcal{GP}(3, \nu, 1)~$ associated with
$~\overline{\mathbb L}(E_{\nu})$, with $~\nu=0, 1, 2$.
Recall that, in most parts of this work,
minimality of the realization is not assumed.

Consider the system matrix
$~L=${\mbox{\tiny$\begin{pmatrix}a&~a&b_1\\ a&~a&b_2\\
b_1&-b_2&0\end{pmatrix}$}} with $~a, b_1, b_2~$ real
parameters. It realizes the function
$~f(s)=\frac{(b_1^2-b_2^2)(s-a)}{s(s-2a)}$.
\vskip 0.2cm

%
%
%
First, for $0>a$, $b_1^2>b_2^2$,
\mbox{$f\in\mathcal{GP}_{\rm min}(3,2,1)\subset\mathcal{P}$},
for $~a>0$,
$~b_2^2>b_1^2$, \mbox{$f\in\mathcal{GP}_{\rm min}(3,0,1)$}
and for $~a(b_2^2-b_1^2)=0$, $~\psi\in\mathcal{GP}(3, 1, 1)$.
\vskip 0.2cm

Next, if $~0>a(b_2^2-b_1^2)~$ $~f~$ is not generalized
positive, but we show that one can always find a static
output feedback of the form \eqref{F} so that
$~f_{\rm cl}\in\mathcal{GP}(3,1,1)$.
Indeed
consider part II
of the claim with $~E_{1,2}$. Condition II(a) is
satisfied. To check condition II(b) note that the
null-space of $~B^*~$ is given by vectors of the form
{\mbox{\tiny$\begin{pmatrix}-b_2\\~b_1\end{pmatrix}$}}.
Thus this condition now reads, $~2a(b_1^2-b_2^2)\geq 0~$
as required.

Indeed, the closed loop system matrix is
$~L_{\rm cl}=${\mbox{\tiny$\begin{pmatrix}
a+kb_1^2&a-kb_1b_2&b_1\\ a+kb_1b_2&a-kb_2^2&b_2\\
b_1&-b_2&0\end{pmatrix}$}} and for $~b_1\not=\pm{b}_2$,
it corresponds to $~\mathcal{GP}_{\rm min}(3,1,1)~$
functions whenever $~\frac{a}{b_2^2-b_1^2}\geq k$.
For instance, for $~a=1$, $b_1=1$, $b_2=1$ and $k=-1$,
one obtains $~L_{\rm cl}=L_{\alpha}~$ from Example
\ref{InertiaExample}.
}
\qed
\end{example}

\section{concluding remarks}
\label{Sec:10}
\setcounter{equation}{0}

As it is often the case, the introduction of a novel concept
opens the door to new research questions.
Concerning the set $~\mathcal{GP}(r,~\nu,~p)~$
we here mention a sample of four problems.

\vskip 0.2cm

\noindent
1. LMI techniques\\
A comprehensive survey of the LMI approach to the PRL appeared in
\cite{GL}. As already mentioned Proposition \ref{PosEquiv1}
suggests that large parts of LMI techniques can be extended to
$~\mathcal{GP}~$ functions. For example, the use of LMI to render
a closed loop function positive, through static output feedback,
was addressed in literature, see e.g. \cite{HSKu,PFA}. This can be
extended in the spirit of Section \ref{Sec:8}. \vskip 0.2cm

\noindent
2. Lyapunov Order\\
Recall that in Examples \ref{InertiaExample} and
\ref{GP(r,nu,p)} we considered the system matrices
$~L_{\gamma}~$ and $~L_{\xi}=L_{\gamma}^{-1}$ and the
respective $~\mathcal{GP}(3, 1, 1)~$ rational functions
\mbox{$\psi_{\gamma}(s)=\frac{1}{s+1}$} and
$~\psi_{\xi}(s)=\frac{s}{s+1}$. Consider now
$~L_x:=\frac{1}{2}(L_{\gamma}+L_{\xi})=${\mbox{\tiny$
\begin{pmatrix}-\frac{1}{2}&0&1\\~0&1&0\\~1&0&\frac{1}{2}
\end{pmatrix}$}} and the associated rational
function $~\psi_x(s)=\frac{s+\frac{5}{2}}{2(s+\frac{1}{2})}~$.
In \cite{CL5} a (partial) Lyapunov order was introduced in
which $~L_{\gamma}\leq L_x$
\begin{footnote}{meaning that whenever
\mbox{$L_{\gamma}\in\overline{\mathbb L}(H)$},~
for some non-singular Hermitian $H$, it implies that for
the same $H$, $~L_x\in\overline{\mathbb L}(H)$.}
\end{footnote}. The Lyapunov order was recently used in
\cite{MS}. It is of interest to find an interpretation of
this partial order in the framework of the rational
functions $~\psi_{\gamma}(s)~$ and $~\psi_x(s)$.
\vskip 0.2cm

\noindent
3. Model order reduction\\
One can exploit the convex structure of all system
matrices associated with the set
\mbox{$\mathcal{GP}(r, \nu, p)$} to try to introduce
a scheme of model order reduction of uncertain systems
in the spirit of \cite[Section 5]{CL2}.
\vskip 0.2cm

\noindent
4. Realization of {\em Even}~ and ~{\em Odd}~ $~\mathcal{GP}~$
functions.\\
Recall that in the scalar case ~{\em Odd}~ functions map
$~i\R~$ to itself while ~{\em Even}~ $~\mathcal{GP}~$
functions map $~i\R~$ to $~\overline{\R}_+$. Both sets were
addressed in \cite{AL2} in the framework of rational functions.
One can study properties of all {\em Even}~ and ~{\em Odd}~
functions within a prescribed set $~\mathcal{GP}(r,~\nu,~p)$.
\vskip 0.2cm

For example, if one considers (for simplicity only real)
realizations of functions in $~\mathcal{GP}(3,~1,~1)$,
the ~{\em Odd}~ and the ~{\em Even}~ cases can be parameterized
by,
\[
\begin{matrix}
L_{\rm odd}=\left(\begin{smallmatrix}0&~a&~b_1\\ a&~0&~b_2\\
b_1&-b_2&~0\end{smallmatrix}\right)&~&~&
L_{\rm even}=\left(\begin{smallmatrix}-a_1&~a_2&~b\\~
a_2&~a_1&~b\\~b&-b&~d\end{smallmatrix}\right)&~
\begin{smallmatrix}a_1+a_2>0\\ b\not=0\\ d\geq 0.
\end{smallmatrix}\\~&~&~&~&~&~&~\\
\psi_{\rm odd}(s)=\frac{(b_1^2-b_2^2)s}{s^2-a^2}&~&~&
\psi_{\rm even}(s)=\frac{b^2(a_1+a_2)}{a_1^2+a_2^2-s^2}+d&~&~
\end{matrix}
\]
It is interesting to note that in the framework of the
associated Lyapunov equation \eqref{Lyap}, in the ~{\em Odd}~
case $~Q=0$, while in
the ~{\em Even}~ case $~Q~$ is diagonal.


\end{document}